\documentclass{amsart}
\usepackage{amssymb, latexsym}

\newcommand{\R}{{\mathbb R}}

\newcommand{\Orderof}{{\mathcal O}}

\newcommand{\vp}{{\varphi}}

\newcommand{\norm}[2]{\|#1\|_{#2}}
\newcommand{\ptl}[2]{\frac{\partial #1}{\partial #2}}
\newcommand{\japanese}[1]{\langle #1\rangle}
\newcommand{\la}{\langle}
\newcommand{\ra}{\rangle}
\newcommand{\les}{\lesssim}

\newcommand{\disp}{\displaystyle}

\theoremstyle{plain}
\newtheorem{theorem}{Theorem}

\newtheorem{proposition}[theorem]{Proposition}
\newtheorem{lemma}[theorem]{Lemma}
\newtheorem{corollary}[theorem]{Corollary}

\theoremstyle{remark}

\begin{document}
\title{Transport in the one-dimensional Schr\"odinger equation}
\author{Michael Goldberg}
\address{Department of Mathematics \\ Johns Hopkins University
\\ 3400 N. Charles  St. \\ Baltimore, MD 21218}
\email{mikeg@math.jhu.edu}

\subjclass[2000]{Primary: 35Q40; Secondary: 34L25}
\keywords{Schr\"odinger equation, dispersive estimates, transport, Jost 
solutions, scattering theory}

\begin{abstract}
We prove a dispersive estimate for the Schr\"odinger equation on the real line,
mapping between weighted $L^p$ spaces with stronger time-decay 
($|t|^{-\frac32}$ versus $|t|^{-\frac12}$) than is possible
on unweighted spaces.  To satisfy this bound, the long-term behavior of
solutions must include transport away from the origin.
Our primary requirements are that
$\langle x\rangle^{3}V$ be integrable and $-\Delta + V$ not have a resonance
at zero energy.  If a resonance is present (for example, in the free case),
similar estimates are valid after projecting away from a rank-one subspace
corresponding to the resonance.
\end{abstract}

\maketitle

In one dimension, the linear propagator of the free Schr\"odinger equation
is given by the explicit convolution
\begin{equation*}
e^{-it\Delta}\psi(x) = \frac{1}{\sqrt{-4\pi i\,t}}\int_\R e^{-i\frac{|x-y|^2}
 {4t}}\psi(y)\, dy.
\end{equation*}
This gives rise immediately to the dispersive estimate
\begin{equation} \label{eq:freedispersive}
\norm{e^{-it\Delta}\psi}{\infty} \le (4\pi |t|)^{-\frac12}\norm{\psi}{1}.
\end{equation}

Such an estimate cannot be true in general for the perturbed operator 
$H = -\Delta + V(x)$.  Even small perturbations of the Laplacian may lead to 
the formation of bound states, i.e.  functions $f_j \in L^2$ 
satisfying $Hf_j = -E_jf_j$.  Bound states with strictly negative energy
are known to possess exponential decay, hence they belong to the entire range
of $L^p(\R)$, $1 \le p \le \infty$.  For each of these bound states $f_j$, 
the associated evolution 
$e^{itH}f_j = e^{-itE_j}f_j$ clearly violates \eqref{eq:freedispersive}.

It is well known \cite{CoddingtonLevinson, ReedSimon} that if $V \in L^1(\R)$
then the pure-point spectrum of $H$ consists of at most countably many
eigenvalues $-E_j < 0$.  The absolutely continuous spectrum of $H$
is the entire positive half-line, and there is no singular continuous spectrum.
Bound states can therefore be removed easily via a spectral projection,
suggesting that one should look instead for dispersive estimates of the form
\begin{equation} \label{eq:dispersive}
\norm{e^{itH}P_{ac}(H)\psi}{\infty} \les |t|^{-\frac12}\norm{\psi}{1}.
\end{equation}

The condition $V \in L^1$ does not always guarantee regularity at the endpoint
of the continuous spectrum.  We say that zero is a resonance of $H$ if there
exists a bounded solution to the equation $Hf = 0$.  Since resonances are
not removed by the spectral projection $P_{ac}(H)$, the validity of 
dispersive estimates invariably depends on whether zero is a resonance of $H$.
Weder~\cite{Weder} and Goldberg-Schlag~\cite{GoldbergSchlag} have shown that
\eqref{eq:dispersive} holds for all potentials with $(1+|x|)^2V \in L^1$, 
and that $(1+|x|)V \in L^1$ suffices provided zero is not a resonance.

The relatively slow time-decay of these estimates (the tail of the function
$t^{-\frac12}$ is not integrable)
makes them unsuitable for many applications.  
We are therefore interested proving a
dispersive estimate which improves the rate of decay by mapping 
between favorably weighted spaces.  Statements of this type appear 
in the work of
Murata \cite{Murata} and Buslaev-Perelman \cite{BuslaevPerelman}, with
weighted $L^2(\R)$ as the underlying space.  A weighted
 $L^1\to L^\infty$ bound was proven recently by Schlag \cite{SchlagSurvey}.
Our first theorem is a refinement of Schlag's result.

\begin{theorem} \label{thm:transport}
Suppose $(1+|x|)^3V \in L^1(\R)$ and zero is not a resonance of $H$.
The continuous part of the Schr\"odinger evolution satisfies the bound
\begin{equation} \label{eq:transport}
\norm{(1+|x|)^{-1} e^{itH}P_{ac}(H)\psi}{\infty} \les |t|^{-\frac32}
\norm{(1+|x|)\psi}{1}.
\end{equation}
\end{theorem}

Recall that $|e^{itH}P_{ac}(H)\psi|(x)$ is always dominated by 
$|t|^{-\frac12}$, by \eqref{eq:dispersive}.  The additional estimate 
\eqref{eq:transport} reduces the bound
even further for all $|x| \ll |t|$.  This suggests that solutions experience
transport away from the origin with nonzero velocity.

The assumption that zero energy is not a resonance is a necessary part
of Theorem~\ref{thm:transport}.  To give an explicit example, consider the
case $V = 0$ with inital data $\psi(x) = e^{-\frac{|x|^2}{2}}$.  For each $t$,
the solution 
$e^{-it\Delta}\psi(x) = (4\pi(1-it))^{-\frac12}e^{\frac{-|x|^2}{2(1-it)}}$
satisfies \eqref{eq:dispersive} but clearly violates \eqref{eq:transport}.
There is a significant degree of structure to a resonance at $\lambda = 0$,
as is seen in the
power-series resolvent expansion of Jensen-Nenciu \cite{JensenNenciu}:
\begin{equation*}
(H - (\lambda + i0))^{-1} = \lambda^{-\frac12}C_{-1} + C_0 + 
\lambda^{\frac12}C_1 + \Orderof(\lambda).
\end{equation*}
Here $C_{-1}$ is a projection onto the subspace spanned by the bounded 
solution of $Hf = 0$, or is vacuous if zero is not a resonance.  
One consequence is that the worst time-decay 
must be confined to a rank-one subspace of functions.  More precisely, in the
one-dimensionsal setting we prove the following:
\begin{theorem}\label{thm:resonance}
Suppose $(1+|x|)^4V \in L^1(\R)$ and there is a nontrivial bounded function 
$f_0$ for which $Hf_0 = 0$, normalized so that $\disp \lim_{x\to\infty} 
(|f_0(x)|^2 + |f_0(-x)|^2) = 2$. 
 Denote by $P_0$ the projection onto the span of $f_0$ given formally by
$P_0\psi = \la \psi,f_0\ra f_0$.

The continuous part of the Schr\"odinger evolution satisfies the bound
\begin{equation} \label{eq:resonance}
\norm{(1+|x|)^{-2} \big(e^{itH}P_{ac}(H)- (-4\pi i\,t)^{-\frac12}P_0\big)
\psi}{\infty} \les  |t|^{-\frac32}\norm{(1+|x|^2)\psi}{1}.
\end{equation}
\end{theorem}

The proof of each theorem relies on a decomposition of the propagator
$e^{itH}P_{ac}(H)$ acoording the the spectral measure of $H$.  Written
this way,
\begin{equation*}
e^{itH}P{ac}(H)\psi = \int_0^\infty e^{it\lambda} E_{ac}(d\lambda) \psi
\,d\lambda
\end{equation*}
where $E_{ac}(d\lambda)$ denotes the absolutely continuous part of the 
spectral measure of $H$.  Since $V$ is assumed to be integrable, it is 
correct to assume that the absolutely continuous spectrum is supported on
the interval $[0,\infty)$.  The Stone formula provides additional information
about the nature of $E_{ac}(d\lambda)$, namely
\begin{equation*}
\la E_{ac}(d\lambda)f,g\ra = \frac1{2\pi i}\la[R_V^+(\lambda) -
 R_V^-(\lambda)]f, g\ra
\end{equation*}
where $R_V^\pm(\lambda) := (-\Delta + V - (\lambda \pm i0))^{-1}$ is the
continuation of the resolvent onto the positive real half-line.
Substituting this into the previous equation yields
\begin{equation*}
\la e^{itH}P_{ac}(H)\psi, \vp\ra = \frac1{2\pi i}\int_0^{\infty} e^{it\lambda}
\la[R_V^+(\lambda) - R_V^-(\lambda)]\psi, \phi\ra\, d\lambda.
\end{equation*}
It is convenient to make the change of variables $\lambda \mapsto \lambda^2$
For the purpose of changing variables inside the resolvent, recall that 
$R_V^+(\lambda)$ is an analytic continuation of the operator-valued
function $(H-z)^{-1}$
from the upper half-plane.  The continuation of $(H-z^2)^{-1}$ is
therefore $(H-(\lambda+i0)^2)^{-1}$, which is identical to 
$R_V^+(\lambda^2)$ along the positive half-line and
$R_V^-(\lambda^2)$ along the negative half-line.  This allows us to
open up the domain of integration to the entire real line:
\begin{equation}\label{eq:integral}
\la e^{itH}P_{ac}(H) \psi, \vp\ra = \frac1{\pi i}\int_{-\infty}^{\infty}
e^{it\lambda^2}\lambda \la R_V^+(\lambda^2)\psi, \vp\ra\, d\lambda.
\end{equation}

For large values of $\lambda$ we will regard $R_V^+(\lambda^2)$ as a 
perturbation of the free resolvent $R_0^+(\lambda^2)$, which can be
expressed explicitly as a convolution.  This part of the argument has
appeared previously in \cite{SchlagSurvey} and requires no further
modification.  For small $\lambda$ we will characterize the resolvent
in terms of the Jost solutions of $H$.  The desired estimates will 
follow from scattering-theory results of Deift and 
Trubowitz~\cite{DeiftTrubowitz}, using similar arguments to those in
Goldberg-Schlag~\cite{GoldbergSchlag}.

To separate the cases of low and high energy, let $\chi$ be a smooth even
cutoff function that is equal to one when $|x| \le \lambda_0$ and is supported
on the interval where $|x| \le 2\lambda_0$.  The value of $\lambda_0$
will be determined later, and depends primarily on the size of $V$.
We will adopt the following notation for discussing polynomially 
weighted $L^p$ spaces.
\begin{align*}
\japanese{x} &:= (1+|x|^2)^{\frac12} \\
\norm{f}{L^p_\sigma} &:= \norm{\japanese{x}^\sigma f}{L^p}
\end{align*}

\section{High Energy Estimates}

Both Theorem~\ref{thm:transport} and Theorem~\ref{thm:resonance} rely on
the same estimate for the high-energy part of the evolution.  This result 
can be found in \cite{SchlagSurvey} but we include it here for the sake of
completeness.
\begin{proposition} \label{prop:highenergy}
Assume that $V \in L^1_1(\R)$ and choose $\lambda_0 \ge \norm{V}{1}$.
The following estimate is valid for all functions $\psi, \vp \in L^1_1(\R)$.
\begin{equation}
\la e^{itH}(1-\chi(\sqrt{H}))\psi, \vp\ra \ \les \ |t|^{-\frac32}\,
\norm{\japanese{x}\psi}{1} \,\norm{\japanese{x}\vp}{1}.
\end{equation}
\end{proposition}

\begin{proof}
By the same spectral argument that led to \eqref{eq:integral}, we are 
estimating here the integral
\begin{equation*}
\frac{1}{\pi i}\int_{-\infty}^\infty e^{it\lambda^2}\lambda (1-\chi(\lambda))
\la R_V^+(\lambda^2)\psi, \vp\ra\, d\lambda.
\end{equation*}
Integrate by parts once to obtain
\begin{equation} \label{eq:highenergy}
\frac{-1}{2\pi t}\int_{-\infty}^\infty e^{it\lambda^2} \frac{d}{d\lambda}
\big[(1 - \chi(\lambda))\la R_V^+(\lambda^2)\psi, \vp\ra\big]\,d\lambda.
\end{equation}

The perturbed resolvent $R_V^+(\lambda^2)$ can be linked to the free
resolvent $R_0^+(\lambda^2)$ via the identity
$R_V^+(\lambda^2) = R_0^+(\lambda^2)(I + VR_0^+(\lambda^2))^{-1}$,
leading to the Born series expansion
\begin{equation*}
R_V^+(\lambda^2) = \sum_{k=0}^\infty R_0^+(\lambda^2)(VR_0^+(\lambda^2))^k.
\end{equation*}

The free resolvent $R_0^+(\lambda^2)$ has an explicit representation as
an integral operator with kernel $K(x,y) = (2i\lambda)^{-1}e^{i\lambda|x-y|}$.
Substituting this into the identity above leads to the
expression
\begin{multline*}
\la R_V^+(\lambda^2) \psi, \vp\ra = \sum_{k=0}^\infty (2i\lambda)^{-(k+1)} 
\int_{\R^{k+2}} e^{i\lambda (\sum_{j=0}^k |x_{j+1}-x_j|)} \\
\psi(x_0) V(x_1)V(x_2)\ldots V(x_k) \vp(x_{k+1})\, dx_0 \ldots dx_{k+1}
\end{multline*}
which is a convergent series provided $2|\lambda| \ge \norm{V}{1}$.
The absence of a boundary term in the integration by parts 
\eqref{eq:highenergy} is justified by a similar argument.

When this is substituted back into the integral \eqref{eq:highenergy}
the differentiation in $\lambda$ leads to two distinct terms.
Up to a constant factor, we have
\begin{align}
\la e^{itH}(1&-\chi(\sqrt{H}))\psi, \vp\ra = \notag \\
&\frac{1}{t}\int_{\R^{k+3}} \sum_{k=0}^\infty \sum_{m=0}^k 
 e^{it\lambda^2}e^{i\lambda (\sum_{j=0}^k |x_{j+1}-x_j|)}
\frac{i(1-\chi(\lambda))}
{(2i\lambda)^{k+1}} |x_{m+1}- x_m| \label{eq:highenergymain}\\
  &\hskip 1.5in \times  \psi(x_0)V(x_1) \ldots V(x_k)\vp(x_{k+1})\, 
  dx_0 \ldots dx_{k+1}\, d\lambda \notag\\
& \hskip -.15in - \frac{1}{t}\int_{\R^{k+3}} \sum_{k=0}^\infty 
e^{it\lambda^2}e^{i\lambda (\sum_{j=0}^k |x_{j+1}-x_j|)} 
\Big[\frac{(k+1)(1-\chi(\lambda))} {\lambda(2i\lambda)^{k+1}} 
 + \frac{\chi'(\lambda)}{(2i\lambda)^{k+1}}\Big] 
  \tag{\ref{eq:highenergymain}a} \label{eq:highenergyother} \\
& \hskip 1.5in \times \psi(x_0)V(x_1) \ldots V(x_k) \vp(x_{k+1})\, 
  dx_0 \ldots dx_{k+1}\, d\lambda.  \notag
\end{align}

In each of the terms we may rearrange the order of integration to
handle the $d\lambda$ integral first (for the $k=0$ term in 
\eqref{eq:highenergymain} this requires restricting to compact support in 
$\lambda$ and taking limits; otherwise it is permitted by Fubini's theorem).
Evaluate this integral using Plancherel's identity:
The Fourier transform of the oscillatory component $e^{it\lambda^2}
e^{i\lambda \sum|x_{j+1}-x_j|}$ is bounded above by $|t|^{-\frac12}$
uniformly in the choice of all $x_j$.  The Fourier transform of each
expression involving the cutoff function (e.g. $\frac{1-\chi(\lambda^2)}
{(2\lambda)^{k+1}}$) is integrable with $L^1(\R)$ norm bounded by
$k (2\lambda_0)^{-(k+1)}$.  This implies that
\begin{align*}
\big|\la e^{itH}(1-\chi&(\sqrt{H})) \psi, \vp\ra\big| \ \les\ 
  |t|^{-\frac32} \sum_{k=0}^\infty k(2\lambda_0)^{-(k+1)} \\
 & \times\int_{\R^{k+2}} \sum_{m=0}^k \la x_{m+1} - x_m\ra\, 
|\psi(x_0)V(x_1)\ldots  V(x_k)\vp(x_{k+1})|\, dx_0 \ldots dx_{k+1}.
\end{align*}
The sum of differences $\la x_{m+1} - x_m\ra$ can be controlled by 
$2\sum_{m=0}^k \la x_m\ra$ using the triangle inequality. The inner integral
is then separable, with the eventual bound
\begin{align*}
\big|\la e^{itH}(1-\chi(\sqrt{H})) \psi, \vp\ra\big| \ &\les\ 
  |t|^{-\frac32} \sum_{k=0}^\infty k^2 (2\lambda_0)^{-k-1}
\norm{V}{1}^{k-1} \norm{\japanese{x}V}{1} \norm{\japanese{x}\psi}{1} 
\norm{\japanese{x}\vp}{1} \\
   &\les \ \lambda_0^{-1}|t|^{-\frac32} \norm{\japanese{x}\psi}{1}
                                        \norm{\japanese{x}\vp}{1}
\end{align*}
provided $\lambda_0 \ge \norm{V}{1}$.
\end{proof}

\section{Low Energy Estimates}
It remains to control the behavior of $e^{itH}\chi(\sqrt{H})\psi$, 
with the result
depending on whether or not $H$ has a resonance at zero.  The Born series
used previously cannot be made to converge, so we rely instead on a 
characterization of the resolvent in terms of Jost solutions.  For each 
$\lambda \in \R$, let $f_\pm(x, \lambda)$ be the unique functions which satisfy
\begin{equation*}
-f''_\pm(x, \lambda) + (V(x) - \lambda^2)f_\pm(x,\lambda) = 0, \qquad
f_\pm(x, \lambda) = e^{\pm i\lambda x}\ {\rm as} \ x \to \pm\infty
\end{equation*}
and $W(\lambda) := W[f_+( \cdot, \lambda), f_-(\cdot,\lambda)]$ be their
Wronskian.  Define also the Wronskian 
$\tilde{W}(\lambda) = W[f_-(\cdot, \lambda), f_+(\cdot, -\lambda)]$. 
The perturbed resolvent $R_V(\lambda^2)$ is an integral operator 
whose kernel is given by
\begin{equation} \label{eq:kernel}
R_V^+(\lambda^2)(x,y) \ = \ \frac{f_+(x, \lambda)f_-(y, \lambda)}{W(\lambda)}
\end{equation}
for all $x \ge y$, and is symmetric for $x < y$.  

Note that $f_\pm(\cdot, -\lambda)$ solve the same second-order 
differential equation as $f_\pm(\cdot, \lambda)$, hence they must be linearly
dependent.  The coefficients in the relation
\begin{equation} \label{eq:intertwining}
f_-(x,\lambda) = \alpha(\lambda)f_+(x, \lambda) 
  + \beta(\lambda)f_+(x,-\lambda)
\end{equation}
are given by 
$\alpha(\lambda) = \frac{\tilde{W}(\lambda)}{-2i\lambda}$ and 
$\beta(\lambda) = \frac{W(\lambda)}{-2i\lambda}$.
These in turn are closely linked to the reflection and transmission
coefficients, namely: $\alpha(\lambda) = \frac{R_1(\lambda)}{T(\lambda)}$
and $\beta(\lambda) = \frac{1}{T(\lambda)}$.  Conjugate symmetry requires
that $\beta(-\lambda) = \overline{\beta(\lambda)}$ and $\alpha(-\lambda)
= \overline{\alpha(\lambda)}$.  Conservation of energy additionally requires
that $|\alpha(\lambda)|^2 + 1 = |\beta(\lambda)|^2$ for every value of 
$\lambda$.  

Since $\beta(\lambda)$ is always positive, $W(\lambda)$ cannot
vanish except possibly when $\lambda = 0$.  The condition $W(0) = 0$ is
satisfied precisely if zero is a resonance; in the generic (non-resonant) case
the values of $W(\lambda)$ are everywhere nonzero.

It is common to rewrite the
Jost solutions as $f_\pm(x, \lambda) = e^{\pm i\lambda x}m_\pm(x,\lambda)$,
where $m_\pm(x,\lambda) \to 1$ as $x \to \pm\infty$.  The relevant properties
of the functions $m_\pm(x,\lambda)$ are summarized below.  See
\cite{DeiftTrubowitz}, Lemma 3 for details.
\begin{lemma} \label{lem:scattering}
Suppose $V \in L^1_\sigma$, $\sigma \ge 1$.  For each $x$ the functions 
$m_\pm(x, \cdot)-1$ belong to the Hardy space $H^{2+}$ of analytic functions
on the upper half-plane.  Consequently, their Fourier transform in the
second variable, denoted by $m_\pm(x, \hat{\rho})$, is supported on the
halfline $\rho \ge 0$.

Define $I(\rho) := \int_{|t| > \rho}|V(t)|\,dt$. The following pointwise
estimates for $m_\pm(x, \hat{\rho})$ are valid over the specified ranges of $x$
and all $\rho > 0$.

\begin{equation} \label{eq:scattering}
{\it If}\ x \ge 0,\ {\it then} \ 
\begin{cases} |m_+(x, \hat{\rho}) - \delta_0(\rho)|
\  \les \ I(\rho) \\
\big|\ptl{}{x}m_+(x, \hat{\rho})\big| \ \les \ I(\rho) + |V(x+\rho)| \\
\big|\ptl{}{\rho}m_+(x, \hat{\rho})\big| \ \les \ I(\rho) + |V(x+\rho)|
\end{cases}
\end{equation}
\begin{equation*}
{\it If}\ x \le 0,\ {\it then} \ 
\begin{cases} |m_-(x, \hat{\rho}) - \delta_0(\rho)|
\  \les \ I(\rho) \\
\big|\ptl{}{x}m_-(x, \hat{\rho})\big| \ \les \ I(\rho) + |V(x-\rho)| \\
\big|\ptl{}{\rho}m_-(x, \hat{\rho})\big| \ \les \ I(\rho) + |V(x-\rho)|
\end{cases}
\end{equation*}
It follows that each of the above functions involving $m_\pm(x, \hat{\rho})$
belongs to $L^1_{\sigma-1}(\R)$, uniformly over all $x$ in the 
appropriate halfline.
Furthermore, the Fourier transform of $\partial_\lambda m_\pm(x, \cdot)$
belongs to $L^1_{\sigma-2}(\R)$.
\end{lemma}
\begin{corollary}
Suppose $V \in L^1_\sigma(\R)$, $\sigma \ge 1$ and let $\tilde{\chi}(\lambda)
= \chi(\frac{\lambda}{4})$.
The functions $\tilde{\chi}(\lambda)W(\lambda)$ and $\tilde{W}(\lambda)$ 
both have Fourier transform (with respect to $\lambda$)
in the space $L^1_{\sigma-1}(\R)$.
\end{corollary}
\begin{proof}
Recall that $f_\pm(x, \lambda) = e^{\pm i\lambda x}m_\pm(x, \lambda)$.
By this definition,
\begin{align*}
\tilde{\chi}(\lambda)W(\lambda) = \tilde{\chi}(\lambda)\big(m_+(0,\lambda)
\partial_x m_-(0,\lambda) - \partial_x &m_+(0,\lambda)m_-(0,\lambda)\big)
\\
&-2i\lambda\tilde{\chi}(\lambda) m_+(0,\lambda)m_-(0,\lambda)
\end{align*}
and
\begin{align*}
\tilde{W}(\lambda) = m_-(0,\lambda)\partial_x
m_+(0, -\lambda) - \partial_x m_-(0,\lambda)m_+(0,-\lambda)
\end{align*}
According to the pointwise bounds in \eqref{eq:scattering}, 
each individual function $m_\pm(0,\pm\lambda)$ has Fourier transform in 
$L^1_{\sigma-1}$, which is an algebra with respect to convolutions.
\end{proof}

\begin{proof}[Proof of Theorem~\textup{\ref{thm:transport}}]
The desired bounds have already been established in the high energy case by 
Proposition~\ref{prop:highenergy}.  The remaining task is to evaluate the part
of the integral not considered in \eqref{eq:highenergy}, namely
\begin{equation*}
 \frac{-1}{2\pi t} \int_{-\infty}^\infty e^{it\lambda^2}
 \frac{d}{d\lambda} \big[\chi(\lambda) \la R_V^+(\lambda^2)\psi, \vp\ra\big]
  \,d\lambda.
\end{equation*}

After applying the formula \eqref{eq:kernel} for the integral kernel of 
$R_V^+(\lambda^2)$ and Plancherel's identity, it suffices to show that
the Fourier transform (in $\lambda$) of
\begin{equation} \label{eq:lowenergy}
\frac{d}{d\lambda} \Big[\chi(\lambda) 
\frac{f_-(x, \lambda)f_+(y, \lambda)}{\tilde{\chi}(\lambda)W(\lambda)} \Big]
\end{equation}
belongs to $L^1(\R)$ with norm bounded by $\japanese{x}\japanese{y}$
for all choices of $x \le y$.  The correct estimate will also hold for
$x > y$ by symmetry of the resolvent.

First consider the case $x \le 0 \le y$.  We are interested in the Fourier
transform of the function
\begin{align*}
i(y-x)\frac{e^{i(y-x)(\cdot)}\chi\,m_-(x,\cdot)m_+(y,\cdot)}{\tilde{\chi}W} &+ 
 \frac{e^{i(y-x)(\cdot)}\partial_\lambda[\chi\, m_-(x,\cdot)
   m_+(y,\cdot)]}{\tilde{\chi}W}  \\
&- \frac{e^{i(y-x)(\cdot)}
    \chi\, m_-(x,\cdot)m_+(y,\cdot) 
     \partial_\lambda[\tilde{\chi}W]}{(\tilde{\chi}W)^2}.
\end{align*}

Lemma~\ref{lem:scattering} ensures that the Fourier transform of each numerator
has $L^1$ norm bounded uniformly in $x \le 0 \le y$.  If zero
is not a resonance, then the Wronskian $W(\lambda)$ is everywhere nonzero.
The Wiener Lemma (see, for example, \cite{Katznelson}, Chapter VIII) then
implies that $\chi(\frac{\lambda}{2})(\tilde{\chi}W)^{-1}$ also has integrable
Fourier transform, making the division possible as well.  Collectively,
the $L^1$ norm of the Fourier transform will be bounded by $|y-x|$ plus
a constant, which in turn is bounded by $\japanese{x}\japanese{y}$.

In the case $0 < x < y$, there is no uniform control over quantities
derived from $m_-(x, \lambda)$.  To avoid this problem, use the 
intertwining coefficients to write
\begin{align*}
f_-(x,\lambda)\  =\ &\alpha(\lambda)f_+(x,\lambda) 
       + \beta(\lambda)f_+(x,-\lambda) \\
= \ &\frac{-1}{4i}\Big[\frac{\tilde{W}(\lambda) + W(\lambda)}{\lambda}
  (f_+(x, \lambda) + f_+(x,-\lambda)) \\
  &\ + (\tilde{W}(\lambda) - W(\lambda))\Big(e^{i\lambda x}
 \frac{m_+(x, \lambda) - m_+(x, 0)}{\lambda} \\
 & \hspace{2.3cm}  + e^{-i\lambda x}
 \frac{m_+(x,0) - m_+(x,-\lambda)}{\lambda}
+ 2i\frac{\sin(\lambda x)}{\lambda}m_+(x,0)\Big)\Big].
\end{align*}
The only functions here of any concern are the expressions with $\lambda$
in the denominator.  Observe that not only is 
$m_+(x, \hat{\rho}) - m_+(x,0)\delta_0(\hat{\rho}) \in L^1_{\sigma-1}(\R)$,
(if one accepts a delta-function at the origin as integrable), but its
integral over the real line is exactly zero.  Because of this, the Fourier 
transform of $\frac{m_+(x, \lambda) - m_+(x, 0)}{\lambda}$ is given by
\begin{equation} \label{eq:cancellation}
\Big[\frac{m_+(x,\cdot) - m_+(x, 0)}{(\,\cdot\,)}\Big]^\wedge(s)
 = -i \int_s^{\infty} m_+(x,\hat{\rho}) \,d\rho
\end{equation}
which belongs to $L^1_{\sigma-2}(\R)$ uniformly in $x \ge 0$.
The term $(m_+(x,0) - m_+(x,-\lambda)/\lambda$ is treated the same way.

An identical argument holds for the fraction $\frac{\tilde{W}(\lambda) + 
W(\lambda)}{\lambda}$ by expanding out each Wronskian according to its
definition.  One can recognize this as a restatement of the well-known
fact about reflection coefficient at zero energy: $R_1(\lambda) = -1$.

To complete the calculations for these terms as in the previous case, one may
need to deal with derivatives such as 
$\frac{d}{d\lambda}\big[\frac{\tilde{W}(\lambda) + W(\lambda)}{\lambda}\big]$.
The Fourier transform of such a function is in $L^1_{\sigma-3}(\R)$, which
is still integrable provided $\sigma \ge 3$.  Depending on where else the
derivative in \eqref{eq:lowenergy} may fall, one obtains norm bounds of size
$\japanese{x} + \japanese{y} + 1$, which is again bounded by $\japanese{x}
\japanese{y}$.

For the term with $\sin(\lambda x)/\lambda$, it is best to go back
to the original integral~\eqref{eq:integral}.
Apply Plancherel's identity to the expression
\begin{equation*}
\int_{-\infty}^\infty e^{it\lambda^2}\sin(\lambda x)  \chi(\lambda)
\frac{((\tilde{W}(\lambda) - W(\lambda))f_+(y,\lambda)}{\tilde{\chi}(\lambda)W(\lambda)} m_+(x,0)\,d\lambda\
\end{equation*}
and observe that the Fourier transform of $e^{it\lambda^2}\sin(\lambda x)$
is a multiple of 
\begin{equation*}
t^{-1/2}\Big(e^{-i\frac{(\rho-x)^2}{4t}} - e^{-i\frac{(\rho+x)^2}{4t}}\Big)
\leq t^{-3/2} |\rho| |x|.
\end{equation*}

The previous estimation of~\eqref{eq:lowenergy} is sufficient to show
that the Fourier transform of 
$\frac{\chi(\lambda) f_+(y, \lambda)}{\tilde{\chi}(\lambda)W(\lambda)}
(\tilde{W}(\lambda) - W(\lambda))$
belongs to $L^1_1(\R)$ with norm controlled by $\la y\ra$.  Thus the size
of this term is not more than $t^{-3/2}|x|\japanese{y}$, as desired.

The case $x \le y < 0$ is handled in an identical by using the 
intertwining relation
$f_+(y, \lambda) = -\overline{\alpha(\lambda)}f_-(y, \lambda) + \beta(\lambda)
f_-(y, -\lambda)$ instead of \eqref{eq:intertwining}.
\end{proof}

\begin{proof}[Proof of Theorem~\textup{\ref{thm:resonance}}]
All of the estimates in Lemma~\ref{lem:scattering} are still valid in the
resonant case.  The one fundamental difference is that $W(\lambda)$ vanishes
when $\lambda = 0$ (and at no other $\lambda \in \R$).  
Consequently, the functions $\alpha(\lambda)$ and $\beta(\lambda) = 
\frac{W(\lambda)}{-2i\lambda}$ are both continuous and real-valued at the
origin.  The Fourier transforms of $\alpha$ and $\tilde{\chi}\beta$ lie in 
the space $L^1_{\sigma-2}(\R)$, and moreover $\beta(\lambda) \not= 0$ over 
the entire real line.

Thanks to the resonance, the integral
\begin{equation*}
\frac{1}{2\pi }\int_{-\infty}^\infty e^{it\lambda^2}\chi(\lambda)
\frac{f_-(x,\lambda)f_+(y,\lambda)}{\tilde{\chi}(\lambda)\beta(\lambda)}
\,d\lambda
\end{equation*}
must have a stationary phase contribution on the order of $|t|^{-1/2}$.
The integrand is sufficently regular that one can isolate the leading term
\begin{align} \label{eq:statphase}
\frac{1}{2\pi}\int_{-\infty}^\infty e^{it\lambda^2}
 \frac{f_-(x, 0)f_+(y, 0)}{\beta(0)}\,d\lambda \ &= \ 
   (-4\pi i\, t)^{-\frac12}\frac{f_-(x,0)f_+(y, 0)}{\beta(0)} \\
 &= \ (-4\pi i\,t)^{-\frac12} f_0(x) f_0(y) \notag
\end{align}
leaving a remainder of order $|t|^{-\frac32}$.  It is clear that
$f_-(\cdot, 0)$ and $f_+(\cdot, 0 )$ are both scalar multiples of $f_0$.
The limiting values of $f_-(x, 0)$ as $x\to \pm \infty$ are 
$\beta(0)+ \alpha(0)$ and $1$, respectively.  This makes
\begin{equation*}
f_-(x, 0) = \sqrt{1+(\beta(0)+\alpha(0))^2}f_0(x).
\end{equation*}
 A similar argument shows that 
$f_+(y,0) = \sqrt{1+(\beta(0)-\alpha(0))^2}f_0(y)$.  The two square roots
have signs in common if $\beta(0) > 0$ and are of opposite sign
if $\beta(0) < 0$.  The last line of \eqref{eq:statphase} is obtained from this
fact and the identity $\alpha^2(0) + 1 = \beta^2(0)$.

The remainder term is given explicitly by
\begin{equation*}
\frac{1}{2\pi} \int_{-\infty}^\infty e^{it\lambda^2} \big(G_{x,y}(\lambda)-
  G_{x,y}(0)\big)\,d\lambda \ = \ \frac{1}{4\pi i\,t} \int_{-\infty}^\infty
e^{it\lambda^2} \frac{d}{d\lambda}
\Big[\frac{G_{x,y}(\lambda) - G_{x,y}(0)}{\lambda}\Big]\,d\lambda
\end{equation*}
 where   
\begin{align*} G_{x,y}(\lambda) &= e^{i\lambda(y-x)}
\chi(\lambda)
\frac{m_-(x,\lambda)m_+(y,\lambda)}{\tilde{\chi}(\lambda)\beta(\lambda)} \\
 &= e^{i\lambda(y+x)}\frac{\chi(\lambda)\alpha(\lambda)m_+(x,\lambda)
m_+(y,\lambda)}{\tilde{\chi}(\lambda)\beta(\lambda)} 
+ e^{i\lambda(y-x)} \chi(\lambda) m_+(x,-\lambda)m_+(y,\lambda).
\end{align*}
One uses the first formula for $G_{x,y}(\lambda)$ in the case $x \le 0 \le y$
and the second formula when $0 < x\le y$.  There is a third formula,
quite similar to the second, which is useful when $x \le y < 0$.

In order to complete the proof it suffices to bound the $L^1$ norm of the
Fourier transform of $\frac{d}{d\lambda}\big[\frac{G_{x,y}(\lambda)-
   G_{x,y}(0)}{\lambda}\big]$ by the quantity $\japanese{x}^2\japanese{y}^2$.
If we are using the second formula for $G_{x,y}$, it is permissible to 
bound each term separately.  All of these estimates are consequences of the 
general rule stated below.

\begin{proposition}
Suppose the Fourier transform of $F(\lambda)$ belongs to $L^1_2(\R)$.
Define 
\begin{equation*}
G(\lambda) = \frac{d}{d\lambda}\Big[\frac{e^{ik\lambda}F(\lambda)-F(0)}
{\lambda}\Big].
\end{equation*}
The Fourier transform of $G$ is integrable, with the bound \ 
$\norm{\widehat{G}}{1} \les \japanese{k}^2\norm{\widehat{F}}{L^1_2}$.
\end{proposition}
Write out $\disp G(\lambda) = 
\frac{d}{d\lambda}\Big[\frac{e^{ik\lambda}-1}{\lambda}\Big] F(\lambda)
+ \frac{e^{ik\lambda}-1}{\lambda} F'(\lambda) + 
\frac{d}{d\lambda}\Big[\frac{F(\lambda)-F(0)}{\lambda}\Big]$.

The Fourier transforms of $\frac{e^{ik\lambda}-1}{\lambda}$ and its derivative
are integrable, with norms proportional to $k$ and $k^2$ respectively.
The Fourier transform of $\frac{F(\lambda)-F(0)}{\lambda}$ belongs to 
$L^1_1(\R)$ (compare to \eqref{eq:cancellation} to see that this is 
controlled by $\norm{\widehat{F}}{L^1_2}$), and that of its derivative is
integrable.  By convolution in $L^1$, each of the three terms above will
yield a bound no greater than $\japanese{k}^2\norm{\widehat{F}}{L^1_2}$,
as desired.
\end{proof}

\bibliographystyle{amsplain}

\end{document}